\documentclass{amsart}

\usepackage{graphicx, amsmath, amsfonts, amssymb, hyperref, amsthm, comment, float, enumerate}
\usepackage{tikz}

\newtheorem{theorem}{Theorem}[section]
\newtheorem{lemma}[theorem]{Lemma}

\theoremstyle{definition}
\newtheorem{definition}[theorem]{Definition}
\newtheorem{example}[theorem]{Example}

\newtheorem{question}[theorem]{Question}

\theoremstyle{remark}
\newtheorem{remark}[theorem]{Remark}

\numberwithin{equation}{section}
\newcommand{\pre}{\textbf{Pre}}
\newcommand{\suf}{\textbf{Suff}}
\newcommand{\sgn}{\textbf{sgn}}
\newcommand{\Next}{\textbf{Next}}
\newcommand{\prev}{\textbf{Prev}}
\newcommand{\wmin}{w_{\text{min}}}

\begin{document}

\title[Concatenations and Poles of Zeta Function]{Concatenations of periodic words in subshifts defined by linear orders and poles of the Artin-Mazur zeta function}


\author[Chenxi Wu]{Chenxi Wu}
\address{Department of Mathematics, UW Madison, Van Vleck Hall, Madison WI 53703}
\email{cwu367@wisc.edu}



\date{\today}

\begin{abstract}
For certain pairs of unimodal maps on the interval with periodic critical orbits, it is known that one can combine them to create another interval map with periodic critical points whose topological entropy is close to one while the poles of Artin-Mazur $\zeta$ function outside the unit circles are close to the other. We provided a formulation and proof of this result in a more general, symbolic dynamical setting, which allow us to generalize this fact to certain families of maps on finite trees.
\end{abstract}

\maketitle

\section{Introduction}

The Milnor-Thurston kneading theory \cite{milnor1988iterated} allows one to establish connections between interval maps and symbolic dynamics, and prior works by Tiozzo \cite{tiozzo2020galois,tiozzo2015topological}, Calegari-Koch-Walker \cite{calegari2017roots}, Lindsey-Tiozzo-Wu \cite{lindsey2021master}, Bray-Davis-Lindsey-Wu \cite{lindsey2023characterization,bray2021shape} and others applied this idea to the study of entropy of certain families of post-critically finite interval maps, as well as core entropies of some superattracting polynomials. In particular, in \cite{bray2021shape,lindsey2021master}, a key result on the dynamics of these families of maps is the following ``persistence theorem'':

\begin{theorem}\cite{bray2021shape, lindsey2021master}
    Let $\mathcal{F}$ be a family of maps, which can be:
    \begin{itemize}
        \item Unimodal maps on intervals with periodic critical orbits.
        \item Maps on Hubbard trees of a superattracting parameter on a principal vein on the Mandelbrot set, with periodic critical orbits.
    \end{itemize}

Every map $f\in F$ admits a finite Markov decomposition. Let $M_f$ be its incidence matrix, $E_f$ the set of eigenvalues of $M_f$, and $h_f$, which is the $\log$ of the leading eigenvalue of $M_f$, the topological entropy.

If $f_1, f_2, f_3\in \mathcal{F}$, $h_{f_1}<h_{f_2}<h_{f_3}$, then for any $\epsilon>0$, any $z$ which is an eigenvalue of $M_{f_1}$ such that $|z|<1$, there is some $f'\in\mathcal{F}$, such that $h_{f_2}<h_{f'}<h_{f_3}$ and $M_{f'}$ has an eigenvalue $z'$ which is $\epsilon$-close to $z$.
\end{theorem}

\begin{figure}[H]
\begin{center}
\includegraphics[scale=0.35]{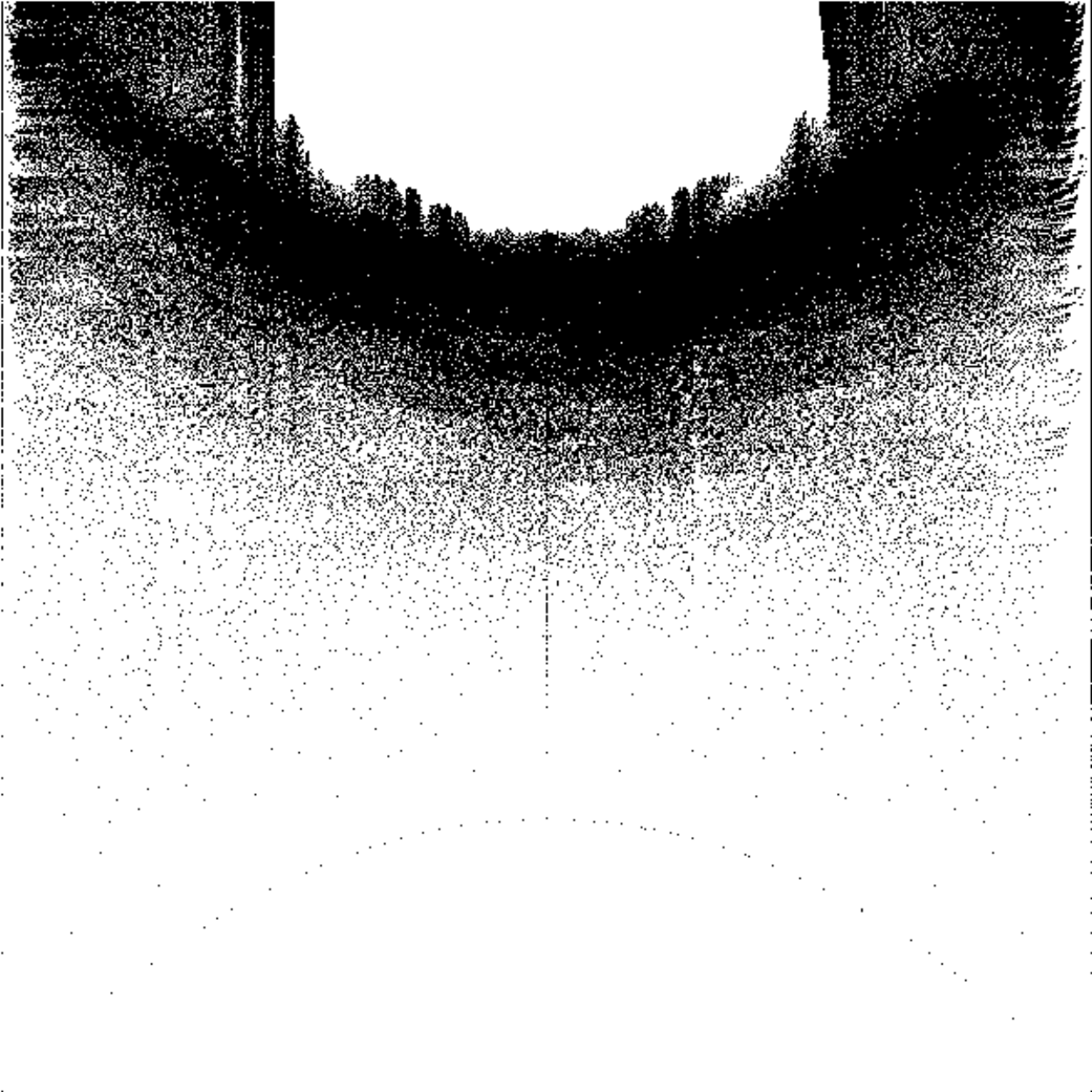}
    \caption{The intersection of Thurston's master teapot with the unit cylinder, plotting the half where the imaginary part is non negative.}
    \label{fig:thurteapot}
\end{center}
\end{figure}

As a consequence, the ``master teapots'' (see \cite{thurston2014entropy}, and Figure \ref{fig:thurteapot}) corresponding to these families of maps, defined as
\[T_\mathcal{F}=\overline{\{(z, \lambda): z\in E_f, \lambda=e^{h_f}, f\in\mathcal{F}\}}\]
when intersecting with the unit cyliner $\{(z, y):|z|<1\}$, consists of vertical line segments that end at the same height. This also provides a connection between horizontal slices of these ``teapots'' with the Mandelbrot set (in the sense of Solomyak \cite{solomyak2005mandelbrot}) of parametrized families of graph directed systems, c.f. \cite{lindsey2023characterization}.

The rough idea behind these kind of ``persistence'' results is the following: one can use Milnor-Thurston kneading theory to relate maps in $\mathcal{F}$ and the forward itinerary of the critical value of these maps, and such itineraries satisfies the condition of {\em admissibility} which means that when deleting any prefix, it becomes no larger than the original itinerary under some linear order. Now for two such maps and corresponding periodic itineraries, such that one has larger entropy, one can find another map in $\mathcal{F}$ with periodic itinerary such that the prefix of the period is close to the itinerary of the map with higher entropy, and the suffix is close to the itinerary of the map with lower entropy. As entropies are controlled by the prefix of the period, while eigenvalues of $M_f$ which are inside the unit circle are controlled by the suffix, one gets the ``persistence'' result.

In this paper, we will examine under exactly what conditions on the linear order we can carry out these kind of  ``concatenation'' argument as in \cite{bray2021shape, lindsey2021master}. From this we will provide a symbolic dynamics formulation of the concept of ``persistence'', which can be generalized to other settings, which will be stated as Theorem \ref{concat2}. In Section 2 we will set the notations, in Section 3-5 we will develop the technical tools of tuning and dominance, in Section 6 we will state and prove the main theorem, and discuss its consequence on Artin-Mazur $\zeta$ functions. Finally, in Section 7 we will provide applications that recovers prior works as special cases as well as applied to families of tree maps. 

In an upcoming paper with Kathryn Lindsey and Ethan Farber, we will use a slight variant of the technique in this paper to show ``persistence'' result for a family of real polynomial maps with 2 critical points.

\section{Notations}

Let $\Gamma=(V, E)$ be a finite directed graph with no parallel edges, that is fully connected and aperiodic. The vertex set is $V=\{0, 1, 2, \dots, N\}$. Let $O: V\rightarrow \{-1, 1\}$ be a function. Let $\Sigma=\{a\in V^{\mathbb{N}}: \text{ for all }n\in \mathbb{N}, (a_n, a_{n+1})\in E\}$ be the set of infinite paths on $\Gamma$, then we have a shift map $\sigma: \Sigma\rightarrow \Sigma$ defined as $(\sigma(a))_n=a_{n+1}$. Here we always assume $0\in\mathbb{N}$.

\begin{definition}
    We define a linear order in $\Sigma$, called the {\em twisted lexicographic order}, as follows: we say $a<b$ if and only if there is some $n\in\mathbb{N}$, such that $a_k=b_k$ for $k<n$, and $\prod_{k=0}^{n-1}O(a_k)(a_n-b_n)<0$. We call an element in $\Sigma$ {\em admissible}, if and only if $a_0=N$ and $\sigma^{n}(a)\leq a$ for all $n\in\mathbb{N}$.
\end{definition}

\begin{remark}
    Let $a\in \Sigma$ be an admissible element, $\Sigma_a=\{b\in \Sigma: b\leq a\}$, then $\sigma(\Sigma_a)\subseteq \Sigma_a$.
\end{remark}

Let $w$ be a finite word with alphabet $V$, i. e. $w\in V^n$ for some $n$. Denote the length of $w$ as $|w|$, and let $\sgn(w)=\prod_k O(w_k)$ be the {\em sign} of $w$. We define a partial order $\leq$ on the set of finite words as follows: if $u, v$ are two finite words, we say $u<v$ iff there is some $0\leq j<\min(|u|, |v|)$, such that $u_k=v_k$ for all $k<j$ and $\prod_{k<j}O(u_k)(u_j-v_j)<0$. The length $k$ prefix and suffix of a word $w$ is denoted by $\pre_k(w)$ and $\suf_k(w)$ respectively.

Some of the following notations follow that of Tiozzo \cite{tiozzo2015topological}. 

\begin{definition}
   Let $w$ be a finite word starting with $N$. We call $w$ {\em periodic} if $w^\infty\in\Sigma$, {\em irreducible} if $w$ can not be written as a concatenation of more than one copies of another word, {\em extremal} if $w$ is periodic and $w^\infty\in \Sigma$ is admissible, {\em admissible} if $w$ is both extremal and has positive $\sgn$, {\em dominant} if  $w$ is periodic, and for any $2\leq k\leq |w|$, $\suf_k(wN)<\pre_k(w)$. Here $wN$ is the concatenation of $w$ with the word $N$.
\end{definition}

\begin{remark}\label{insig}
    A finite word $w$ satisfies $w^\infty\in \Sigma$ iff it satisfies the following two conditions:
    \begin{enumerate}[(1)]
        \item For any $0\leq i\leq |w|-2$, $(w_i, w_{i+1})\in E$.
        \item $(w_{|w|-1}, w_0)\in E$.
    \end{enumerate}
\end{remark}

\begin{remark}
A finite word $w$ starting with $N$ is extremal iff $w^\infty\in \Sigma$, and for any $k<|w|$, $\suf_{|w|-k}(w)\pre_k(w)\leq w$. Also, if $w$ is dominant then it is both extremal and irreducible.
\end{remark}

Let $a, b\in \Sigma$, $a<b$, the {\bf open interval} $(a, b)$ is the elements in $\Sigma$ that lies between $a$ and $b$, and the closed interval is the open interval union with $\{a, b\}$, in other words, $[a, b]=\{c\in \Sigma: a\leq c\leq b\}$. Let $I_k$ be the set of elements in $\Sigma$ that starts with $k$.

\begin{remark}\label{markov}
If $w$ is extremal, then $\Sigma_{w^\infty}$ admits a finite Markov decomposition as follows: Let $a_1>a_2>\dots >a_m$ be distinct elements in the set $\{\sigma^k(w^\infty): k\in\mathbb{N}\}$. Let $a_{m+1}$ be the smallest element in $\Sigma$. Now it is easy to see that $\mathcal{P}=\{[a_{i+1}, a_i]\cap I_j: 1\leq j\leq N, 1\leq i\leq m+1, [a_{i+1}, a_i]\cap I_j\not=\emptyset\}$ form a Markov decomposition of $\Sigma_{w^\infty}$. Hence the Artin-Mazur $\zeta$-function of $\sigma|_{\Sigma_{w^\infty}}$ is rational with integer coefficients, and the topological entropy is the logrithm of an algebraic integer. Furthermore, because if $w_1^\infty<w_2^\infty$ then $\Sigma_{w_1^\infty}\subseteq \Sigma_{w_2^\infty}$, the topological entropy of $\Sigma_{w_2^\infty}$ is no less than the topological entropy of $\Sigma_{w_1^\infty}$. In other words, the topological entropy of $\Sigma_{w^\infty}$ is monotone with respect to $w^\infty$ under the linear order $<$.
\end{remark}

\section{Concatenation of Dominant and Admissible Words}

The purpose of this section is to prove a generalization of Proposition 5.1 of \cite{bray2021shape}, which says that one can create admissible words by concatenating dominant and admissible words. The proof follows closely to that of \cite{bray2021shape} as well.

Firstly some elementary lemmas, which are generalizations of arguments in Carminati-Tiozzo \cite{carminati2012canonical}. See also \cite{tiozzo2015topological,bray2021shape}.

\begin{lemma}\label{flip}
    Let $a, b$ be two periodic words, in other words, $a^\infty, b^\infty\in \Sigma$, and both start at $N$. Then
    \begin{enumerate}[(1)]
        \item $(ab)^\infty$ and $(ba)^\infty$ are both in $\Sigma$.
        \item $a^\infty>b^\infty$ iff $ab>ba$ 
    \end{enumerate}
\end{lemma}

\begin{proof}
Part (1) follows from the observation in Remark \ref{insig}.

    Suppose $ab=ba$, then $a$ and $b$ must be common multiples of the same word, which implies $a^\infty=b^\infty$. Similarly $a^\infty=b^\infty$ also implies $ab=ba$. Now to show (2) one only need to show that it is not possible to have both $a^\infty>b^\infty$ and $ab<ba$. Suppose $a^\infty>b^\infty$ and $ab<ba$. We shall now look at where $ab$ and $ba$ first differ. Suppose they differ for the first time in the $k$-th letter, i.e. $\pre_{k-1}(ab)=\pre_{k-1}(ba)$, and $\pre_k(ab)\not=\pre_k(ba)$.

    \begin{itemize}
        \item Case 1: If $k\leq\min(|a|, |b|)$, this implies $a<b$, hence $a^\infty<b^\infty$, a contradiction.
        \item Case 2: If $|a|=|b|$, then it is easy to see that $k\leq |a|=|b|$, which reduces to case 1.
        \item Case 3: If $|a|<|b|$, and $k>|a|$. This means that $b=a^jc$ for some word $c$ which does not have $a$ as a prefix, and $j>0$ a positive integer. Hence $a^{j+1}c<a^jca$. Now if $k\leq (j+1)|a|$, we can similarly get $a^\infty<b^\infty$. If $k>(j+1)|a|$, then $c$ must be a prefix of $a$, hence $a^{j+1}c$ is a prefix of $a^\infty$ while $a^jca$ is a prefix of $b^\infty$, this also implies that $a^\infty<b^\infty$.
        \item Case 4: If $|b|<|a|$, and $k>|b|$. This can be shown to be impossible via the same argument as in Case 3.
    \end{itemize}
\end{proof}

An immediate consequence of Lemma \ref{flip} is the following:

\begin{lemma}\label{flip2}
    If $a^\infty$ and $b^\infty$ are both in $\Sigma$ and both starts with $N$, then $a^\infty>b^\infty$ iff $ab^\infty>b^\infty$ iff $ba^\infty<a^\infty$. 
\end{lemma}

\begin{proof}
    Apply Lemma \ref{flip} to $a$ and $b^n$ or $a^n$ and $b$, and set $n\rightarrow\infty$.
\end{proof}

\begin{lemma}\label{sand}
    Let $w$ be an extremal and irreducible word, $0<k<|w|$. If $\pre_k(w)=\suf_k(w)$, then $\sgn(\pre_k(w))=-1$.
\end{lemma}

\begin{proof}
    Let $v=\pre_k(w)$, then $w$ can be written as $w=vx=yv$, and extremality and irreducibility of $w$ implies that $vx>xv$, $yv>vy$.

    Now assume that $\sgn(v)=1$, then
    $vx=yv>vy$
    hence $x>y$, and $yv=vx>xv$ hence $y\geq x$, a contradiction. 
\end{proof}

Now we state and prove the main theorem for this Section:

\begin{theorem}\label{concat1}
    Let $w$ be a dominant word, $v$ an admissible and irreducible word, $v^\infty<w^\infty$, and $n$ a positive integer such that $n|w|<2|v|$, then $w^nv^n$ is also admissible.
\end{theorem}

\begin{proof}
By Remark \ref{insig}, $w^nv^n$ is periodic, and it is evident that it has positive $\sgn$. Now we only need to show that $\sigma^k((w^nv^n)^\infty)\leq (w^nv^n)^\infty$, for $1\leq k\leq n|w|+n|v|-1$.

\begin{enumerate}[(1)]
\item If $k=t|w|$ for $1\leq t\leq n$, remove the common prefix $w^{n-t}$, then compare the $n|w|+n|v|$ prefixes by applying the second part of Lemma \ref{flip}.
\item If $k=n|w|+t|v|$ for $0\leq t\leq n-1$, compare the $(n-t)|v|+n|w|$ prefixes by applying the second part of Lemma \ref{flip}.
\item If $k=t|w|+r$, where $0\leq t\leq n-1$, $0<r<|w|$, compare the $|w|+1-r$ prefixes, by making use of the assumption that $w$ is dominant.
\item If $k=n|w|+t|v|+r$ where $0\leq t<n-1$, $0<r<|v|$. 
\begin{enumerate}[(a)]
\item If the $|v|-r$ prefixes are different, 
\[\pre_{|v|-r}((w^nv^n)^\infty)=\pre_{|v|-r}w^\infty\]
\[\pre_{|v|-r}(\sigma^k(w^nv^n)^\infty)=\suf_{|v|-r}v\leq \pre_{|v|-r}v\leq \pre_{|v|-r}w^\infty\]
The first $\leq$ due to $v$ being admissible, and the second $\leq$ due to $v^\infty<w^\infty$. Hence $\sigma^k(w^nv^n)^\infty<(w^nv^n)^\infty$.
\item If the $|v|-r$ prefixes are the same, by Lemma \ref{sand}, this common $|v|-r$ prefix has negative cumulative sign. Deleting this common prefix, we reduce it to case (1) or case (3).
\end{enumerate}
\end{enumerate}
\end{proof}

\section{Tuning and Renormalizability}

\begin{definition} Let $W_n$ be the set of periodic words of length $n$. Then the $\leq$ defined in Section 2 gives a linear order on the finite set $W_n$. For any $w\in W_n$, let $\Next(w)$ be the element in $W_n$ which is right above $w$ under $\leq$, and $\prev(w)$ the element that is right below $w$.
\end{definition}

\begin{definition}
    We say a pair of elements $a, a'$ in $W_n$ form a {\em tuning pair}, if:
    \begin{itemize}
        \item $a$ is admissible, $a'$ is extremal and $\sgn(a')=-1$.
        \item $a'=\Next(a)$.
    \end{itemize}
\end{definition}

\begin{definition}
    We call the tuple $(\Gamma=(V, E), O)$ to be {\em tunable} if
    \begin{enumerate}
        \item There is a periodic word $\wmin$, with only a single $N$, and $\wmin N$ is the minimal element in $W_{|\wmin|+1}$.
        \item If $a\in W_n$ is extremal and irreducible, and $\sgn(a)=-1$, then either $a=\wmin^k$ for some odd number $k$, or $a$ lies in a tuning pair.
    \end{enumerate}
\end{definition}

\begin{example}\label{unim}
    If 
    \[\Gamma_2=(\{0, 1\}, \{(0, 0), (0, 1), (1, 0), (1, 1)\})\]
    $O_2$ sends $0$ to $1$ and $1$ to $-1$. This is the symbolic dynamics corresponding to unimodal interval maps \cite{tiozzo2020galois,bray2021shape,lindsey2023characterization}.

    Now we show that $(\Gamma_2, O_2)$ is tunable. Firstly it is easy to see that $\wmin=1$ satisfies the first condition. To check the second condition, suppose $a\not=1^n$ and $a$ is extremal, irreducible and $\sgn(a)=-1$. Let $b$ be $a$ with the last letter removed, then $\sgn(b)=-\sgn(a_{n-1})$, hence 
    \[a'=\begin{cases}b0 & a_{n-1}=1\\ b1 & a_{n-1}=1\end{cases}\]
    equals $\prev(a)$ and $\sgn(a')=1$. Now we only need to show that $a'$ is also extremal. Suppose otherwise, there are $u$, $v$ such that $a'=uv$, and $vu>uv$, hence $vu\geq a$,
    \[v\geq \pre_{|v|}(a')=\pre_{|v|}(a)\geq \suf_{|v|}(a)\]
    and $v\not=\suf_{|v|}(a)$ as they differ at the last digit, hence $v>\suf_{|v|}(a)$. this, together with the fact that $uv=a'<a=u\suf_{|v|}(a)$, imply that $\sgn(u)=-1$, hence $\sgn(v)=-1$, and $\sgn(\suf_{|v|}(a))=-1$. Because $v$ and $\suf_{|v|}(a)$ differ only at the last letter, there can be no word in between, hence $\pre_{|v|}(a')=\pre_{|v|}(a)$ must equals either $v$ or $\suf_{|v|}(a)$. The latter case contradicts with Lemma \ref{sand}. In the former case, deleting the common $|v|$-prefix of $vu$ and $a$ one gets $u\leq \suf_{|u|}(a)$, which, due to the extremality of $a$, can only mean that $u=\suf_{|u|}(a)$, i.e. $vu=a$. Yet $\sgn(a)=-1$, $\sgn(vu)=\sgn(uv)=\sgn(a)=1$, a contradiction. 
\end{example}

\begin{remark}\label{fliplastdig}
    The argument above means that to check a pair is tunable, one can first check part 1, and next show that for any $w\in W_n$, $\prev(w)$ is $w$ with one letter changed, where the letter is either the last one or is followed by a suffix that is shared by all periodic words. 
\end{remark}

From now on we always assume that our $(\Gamma, O)$ is tunable.

\begin{definition}
    Let $(w, w')$ be a tuning pair, $v$ an admissible word in the setting of $(\Gamma_2, O_2)$ as in Example \ref{unim}. Then by the {\em tuning} of $w$ by $v$ we mean the word
    \[w*v=w_0w_1\dots w_{|v|-1}\]
    such that
    \[w_i=\begin{cases} w & v_i=0\\ w' & v_i=1\end{cases}\]
    If an admissible word can not be written as a tuning we call it {\em nonrenormalizable}. Similarly one can define tuning of admissible infinite sequences, where $v$ is an admissible infinite sequence in the setting of $(\Gamma_2, O_2)$.
\end{definition}

\begin{lemma}\label{minnonnor}
    If $w$ is extremal but not admissible (i.e. $w^\infty$ admissible and $\sgn(w)=-1$), and $w\not=\wmin^n$, then $w$ can be written as a tuning, and there are $w'$, $v$ such that $w=w'*v$, and $w'$ is nonrenormalizable.
\end{lemma}

\begin{proof}
    $w$ can always be written as $w=w_1^k$, where $k$ is an odd integer, $w_1$ is irreducible, extremal and $\sgn(w_1)=-1$. Hence tunability of $(\Gamma, O)$ implies that $(\prev(w_1), w_1)$ is a tuning pair, and $w=\prev(w_1)*1^k$.
    
    Consider the non-empty set $\{(w', v): w=w'*v\}$, pick $(w'_0, v_0)$ such that $w'_0$ is the shortest. Suppose $w'_0=w''*v''$, where $|v''|>1$, then by the definition of the order and the fact that $\Next(v_0)$ is always $v_0$ with the last digit flipped (see Example \ref{unim}), $\Next(w'_0)=w''*\Next(v'')$, and $w=w''*(v''*v_0)$, a contradiction.
\end{proof}

\begin{remark}\label{rmk:irred}
It is easy to see that $v<v'$ iff $w*v<w*v'$. Furthermore, if $w$ is admissible, $a$ is nonrenormalizable, and $w=a*v$, then the intervals in Example \ref{markov} that are of the form $[\sigma^{i}(w^\infty), \sigma^{i+k|a|}(w^\infty)]$ where $k\in \mathbb{Z}$ must be sent to one another, in other words the Markov decomposition incidence matrix can not be irreducible.
\end{remark}

\begin{definition}\label{dfn:uniqmintun}
    We say that the tuple $(\Gamma=(V, E), O)$ has {\em unique maximal tuning}, if the $w'$ in Lemma \ref{minnonnor} is always {\em minimal}, in the sense that if there is some $w''$, $v''$ such that $w=w''*v''$, then the length of $w''$ is a multiple of the length of $w'$.
\end{definition}

\begin{remark}
    One can show that tuples that satisfies Remark \ref{fliplastdig} also has unique maximal tuning.
\end{remark}

From now on we always assume that our $(\Gamma, O)$ has unique maximal tuning.

\section{Construction of Dominant Words}

The main goal of this Section is to prove the following, which is a generalization of Proposition 10.5 of \cite{tiozzo2015topological}:

\begin{theorem}\label{dendom}
    Suppose $w$ is an admissible, irreducible and non-renormalizable word, then either $w$ is dominant, or for any $n$, there is some word $b$ such that $w^nb$ is dominant.
\end{theorem}

Following the notation in \cite{tiozzo2015topological}, for any admissible word $w$, let $PS(w)$ be the set of common prefixes and suffixes of $w$, i.e. 
\[PS(w)=\{v: 0<|v|<|w|, w=yv=vz, z_0=N\}\]
And 
\[RS(w)=\{z: w=vz, v\in PS(w)\}\]
It is evident that elements of $PS(w)$ and $RS(w)$ are all periodic, and an admissible word $w$ is dominant iff $PS(w)=\emptyset$.

Firstly we need the following four lemmas:

\begin{lemma}\label{denselem1}
    Let $w$ be an admissible and nonrenormalizable word which is not dominant. Let $z_0$ be the element of $RS(w)$ that maximizes $z_0^\infty$. Let $v$ be the word such that $w=vz_0$, then $v$ is extremal.
\end{lemma}

\begin{proof}
    Consider any decomposition of $v$ into $v=ab$, we only need to show $ba\leq ab$. By definition of $RS(w)$ and $PS(w)$, $v$ is a suffix of $w$ as well, hence $\sigma^{|w|-|b|}(w^\infty)=bw^\infty<w^\infty$. Comparing the $|v|$-prefix of both sides we see that $ba\leq ab$.
\end{proof}

\begin{lemma}\label{seqtune}
    Let $a$ be a admissible word such that $(a, a')$ form a tuning pair, $\beta$ an admissible infinite sequence such that $a*0^\infty\leq \beta\leq a*10^\infty$, then there is an admissible infinite sequence $s$ in the setting of $(\Gamma_2, O_2)$ such that $\beta=a*s$.
\end{lemma}

\begin{proof}
    We only need to show that each section of $\beta$ from the $k|a|+1$-th place to the $(k+1)|a|$ place, where $k$ is arbitrary natural number, is either $a$ or $a'$.

    By admissibility of $\beta$, this section is no more than the $|a|$-prefix of $\beta$, which is no more than $a'$. Now the prior sections, by inductive hypothesis, are all either $a$ or $a'$. If before this section there is some $a'a^j$ ($j$ may be $0$), do $\sigma^{(k-1-j)|a|}$ and compare the $(j+2)|a|$ prefix with $\beta$ (which is no more than $a*10$ by assumption), we see that this section can not be less than $a$. If the previous sections are all $a$, by $\beta\geq a^\infty$ we see that this section is no less than $a$ as well.
\end{proof}

\begin{lemma}\label{denselem2}
    Let $w$ be an admissible and nonrenormalizable word which is not dominant. Let $z_0$ and $v$ be defined as in Lemma \ref{denselem1}, and $v=\alpha * t$ for some nonrenormalizable word $\alpha$. Then $\alpha^\infty>z_0^\infty$.
\end{lemma}
\begin{proof}
        Let $\alpha' = \Next(\alpha)$. Because $10^\infty$ is the largest $(\Gamma_2, O_2)$-sequence, while any $(\Gamma_2, O_2)$-sequence is larger than $0^\infty$,
    \[\alpha'\alpha^\infty=\alpha*10^\infty>\alpha*t0^\infty=v\alpha^\infty\]
    \[v^\infty=\alpha*t^\infty>\alpha*0^\infty=\alpha^\infty\]
    Now let $y$ be a word such that $w=yv$, then $y^\infty>v^\infty=(v^n)^\infty$ for any positive integer $n$, hence 
    $yv^n>v^ny$ for any $n$. Let $n$ goes to infinity we have $yv^\infty=wv^\infty>v^\infty$
    which implies that 
    $w^\infty>v^\infty>\alpha^\infty$.

    Now we proceed by looking at various cases:
        \begin{itemize}
            \item Case 1: $(z_0v)^\infty\geq \alpha^\infty$. We have \[w^\infty=v(z_0v)^\infty\leq v\alpha^\infty<\alpha'\alpha^\infty\]
            which implies that $w=\alpha*t'$ for some $(\Gamma_2, O_2)$-word $t'$ by Lemma \ref{seqtune}. Apply Definition \ref{dfn:uniqmintun} we see that this contradicts with the assumption that $w$ is non renormalizable.
            \item Case 2: $(z_0v)^\infty<\alpha^\infty$. There are two subcases:
            \begin{itemize}
                \item Case 2(a): The $|z_0|$-prefix of $\alpha^\infty$ is strictly larger than $z_0$: by looking at the $|z_0|$-prefix we have $\alpha^\infty>z_0^\infty$.
                \item Case 2(b): The $|z_0|$-prefix of $\alpha^\infty$ equals $z_0$. This means that $z_0$ is of the form $z_0=\alpha^kp$ where $p$ is a proper prefix of $\alpha$. By Lemma \ref{sand}, $\sgn(z_0)=\sgn(v)=-1$, so $\sgn(p)=-1$. Now delete from $\alpha^\infty$ and $z_0^\infty$ the common prefix $z_0$, by admissibility and irreducibility of $\alpha$ we see that $\alpha^\infty>z_0^\infty$.
            \end{itemize}
        \end{itemize}
    \end{proof}

\begin{lemma}\label{denselem3}
Let $w$ be an admissible and nonrenormalizable word which is not dominant. Let $u$ be a dominant word, such that
\begin{itemize}
    \item $w>u$.
    \item For any $z\in RS(w)$, $u^\infty>z^\infty$.
\end{itemize}
Then for sufficiently large $n$, $w^nu^n$ is dominant.
    \end{lemma}
    
    \begin{proof}
We need to show that for any $k\geq 1$, the $k+1$-prefix of $w^nu^n$ is strictly larger than the $k$ suffix of $w^nu^n$ concatenated with $N$.

\begin{itemize}
    \item Case 1: $k\leq n|u|$: the $k$ suffix concatenated with $N$ either start with a proper suffix of $u$ followed by $N$, which due to the dominance of $u$, is strictly smaller than the corresponding prefix of $u$, which in turn is no larger than the corresponding prefix of $w$; or it may start with a $u$, which is smaller than $w$ by assumption.
    \item Case 2: $k\geq n|u|+|w|$: the $k$-suffix either starts with a shift of $w$, or is a (less than $n$) multiple of $w$ followed by $u^n$. In the former case, admissibility and irreducibility of $w$ implies that this suffix is strictly smaller than the corresponding prefix. In the latter case, the suffix is strictly smaller than the corresponding prefix because $u<w$.
    \item Case 3: $n|u|<k<n|u|+|w|$. In this case the $k$-suffix is of the form $du^n$ where $d$ is a proper suffix of $w$. 
    \begin{itemize}
        \item Case 3a: If $dN$ is not also a prefix of $w$, $dN$ must be strictly smaller than the corresponding prefix of $w$ due to admissibility of $w$, hence this $k$-suffix must be strictly smaller than the corresponding prefix.
        \item Case 3b: $dN$ is also a prefix of $w$. Let $w=de$, then $e\in RS(w)$ by definition. Because $w^k=(de)^k$ is admissible, we have $(w^{k-1}d)^\infty>e^\infty$. Let $k$ be sufficiently large we have $w^\infty\geq e^\infty$. By Lemma \ref{sand}, $sgn(e)=sgn(d)=-1$, hence $ew^\infty\leq e^\infty<u^\infty$, $dew^\infty=w^\infty>du^\infty$, so when $n$ is large enough, $du^nN$ is strictly smaller than the corresponding prefix in $w^\infty$. 
    \end{itemize}
    \end{itemize}
\end{proof}

Now we are ready to prove the main theorem:

\begin{proof}[Proof of Theorem \ref{dendom}]
    Suppose the theorem is not true, let $w$ be the shortest word which is admissible, nonrenormalizable, non-dominant, and there exists $n \in \mathbb{N}$ such that $w^nr$ is not dominant for any $r$.

    Let $z_0=\arg\max_{z\in RS(w)}z^\infty$,  and let $v$ be a word such that $w=vz_0$. By Lemma \ref{sand}, $\sgn(v)=-1$, and by Lemma \ref{denselem1}, $v$ is extremal. Hence by Lemma \ref{minnonnor}, either of the two cases below are true: 
    \begin{enumerate}[i)]
        \item $v=\alpha * t$ for some nonrenormalizable word $\alpha$.
        \item $v=\wmin^k$, $k$ is odd, $\sgn(\wmin)=-1$.
    \end{enumerate}
    
    In case i), by Lemma \ref{denselem2}, $\alpha^\infty>z_0^\infty$. Because $|\alpha|<|w|$ and  $w$ is the shortest word where the claim does hold hold, either $\alpha$ itself is dominant, or for any integer $n$ there is some word $r'$ such that $\alpha^nr'$ is dominant. We now have a dominant word $u$, which equals $\alpha$ if $\alpha$ is dominant itself, and $\alpha^nr'$ for $n$ sufficiently large if otherwise, such that $w>u$, $u^\infty>z^\infty$ for any $z\in RS(w)$. Here $u^\infty>z^\infty$ is because $u^\infty$ shares a long common prefix with $a^\infty$, which is greater than $z_0^\infty$. Hence by Lemma \ref{denselem3}, when $M\gg0$, $w^Mu^M$ is dominant. So let $r=w^{M-n}u^M$, we get a contradiction.

    In case ii), because $w^\infty\geq \sigma^{|\wmin|}(w^\infty)$, $w$ must equals $\wmin^k$, which contradicts with the assumption of $w$ being irreducible. 
\end{proof}
    
\section{Topological entropy and Poles of the Artin-Mazur zeta function}

Combining the results from the previous 3 sections, we get the following:

\begin{theorem}\label{concat2}
    Suppose $(\Gamma, O)$ is tunable and has unique maximal tuning. For any nonrenormalizable word $a$, any admissible word $b$, any positive integer $n$, there is a word $c$ such that $a^ncb^n$ is admissible.
\end{theorem}

\begin{proof}
    This follows immediately from Theorem \ref{concat1} and Theorem \ref{dendom}.
\end{proof}

Let $\Gamma$ be a fully connected and aperiodic directed graph, with vertices labeled $0, \dots, N$, and let $w$ be an admissible word. Then by Remark \ref{markov}, the subshift $(\Sigma_{w^\infty}, \sigma)$ admits a Markov decomposition, and the poles of the Artin-Mazur $\zeta$ function are reciprocals of the eigenvalues of the Markov decomposition incidence matrix. These eigenvalues can also be written as roots of a polynomial as below:

By Gaussian elimination, there are polynomials $p_0(\lambda), \dots, p_N(\lambda)$, such that for each $0\leq i\leq N-1$, 
\[\lambda p_i(\lambda)=\sum_{(i, j)\in E(\Gamma)}p_j(\lambda)\]
and $\gcd(p_0, \dots, p_N)=1$. Now for any edge $(i, j)$ in $\Gamma$, define
\[F_{i,j}(x)=\begin{cases} \lambda x-\sum_{(i, k)\in E(\Gamma), k<j}p_k(\lambda) & O(i)=1\\p_j(\lambda)-\lambda(p_i(\lambda)-x)+\sum_{(i, k)\in E(\Gamma), k>j}p_k(\lambda) & O(i)=-1\end{cases}\]

And let 
\[F_w(\lambda)=F_{w_{|w|-1}, N}\circ \dots \circ F_{w_0, w_1}(p_N(\lambda))-p_N(\lambda)\]
With exactly the same argument as in Theorem 1.2 of \cite{lindsey2021master}, one can show that the roots of $F_w$ that are off unit circle are exactly the eigenvalues of the Markov decomposition matrix $M$ which are off the unit circle, and furthermore by construction the coefficients of $F_w$ are uniformly bounded with the bound only dependent on $(\Gamma, O)$, while the lowest non-zero degree is bounded as well. Hence, by Rouche's theorem, the eigenvalues outside of the unit circle are controlled by prefix of $w$ while the eigenvalues inside the unit circle are controlled by suffix of $w$. This observation, together with Theorem \ref{concat2} and Remark \ref{rmk:irred}, implies that:

\begin{theorem}\label{concat3}
    Let $(\Gamma, O)$ be tunable and has unique maximal tuning. Let $w$ be an admissible word such that the Markov decomposition matrix of $\Sigma_{w^\infty}$ defined in Remark \ref{markov} is irreducible, and $v$ an admissible word such that $v^\infty<w^\infty$, then, for any $\epsilon>0$, any pole $z$ of the Artin-Mazur $\zeta$ function of $\Sigma_{v^\infty}$ outside the unit circle, there is some admissible word $w'$ such that the entropy of $\Sigma_{w'^\infty}$ is within $\epsilon$ to the entropy of $\Sigma_{w^\infty}$, and a pole of the Artin-Mazur $\zeta$ function of $\Sigma_{w'^\infty}$ is within $\epsilon$ distance to $z$.
\end{theorem}\qed

\begin{remark}
    Similar to \cite{lindsey2023characterization}, Theorem \ref{concat3} also implies that when $w$ is a nonrenormalizable admissible word, the suffix of all possible admissible words $w'$ such that $w'^\infty\leq w^\infty$, when reversed order and taken closure, is a subshift of finite type. Furthermore, the closure of the set of all Markov matrix eigenvalues for such $\Sigma_{w'^\infty}$ which are inside the unit circle, can be written as the limit set of graph directed system, hence can be studied via the technique of Solomyak \cite{solomyak2005mandelbrot}.
\end{remark}

\section{Applications to Dynamics and further questions}

\begin{example}
    By Milnor-Thurston kneading theory \cite{milnor1988iterated}, and Remark \ref{unim}, Theorem \ref{concat3} implies Theorem 1 in \cite{bray2021shape} for $\lambda>\sqrt{2}$. By Branner-Hubbard surgery \cite{branner1988surgery, lindsey2021master, riedl2001arcs} one can obtain the ``persistence'' result in \cite{lindsey2021master} as well.
\end{example}

\begin{example}
    Let $\Gamma$ be a graph with vertices $\{0, 1, \dots, n\}$, $O(k)=(-1)^k$, and there is a directed edge from any vertex to any vertex. This $(\Gamma, O)$ is tunable with $\wmin=n$, and $\prev(w)$ is $w$ with the last letter $i$ replaced by either $i+1$ or $i-1$. This corresponds to continuous interval maps that have $n$ critical points and $2$ critical values, where the smaller critical value is also a fixed point.
\end{example}

\begin{example}
     Let $\Gamma$ be a graph with vertices $\{0, 1, 2, 3\}$, and there are edges from $i$ to $j$ except for $i=1$ or $2$ and $j=3$, let $O(k)=(-1)^k$. Then we can check that $(\Gamma, O)$ is tunable:
     \begin{enumerate}
         \item $\wmin=3$.
         \item If $w$ is an extremal word that is not admissible, $\prev(w)$ is $w$ with the last $0$ replaced with $3$ or vice versa. Now apply Remark \ref{fliplastdig}.
     \end{enumerate}
\end{example}

\begin{figure}[H]
\begin{center}
\includegraphics[scale=0.35]{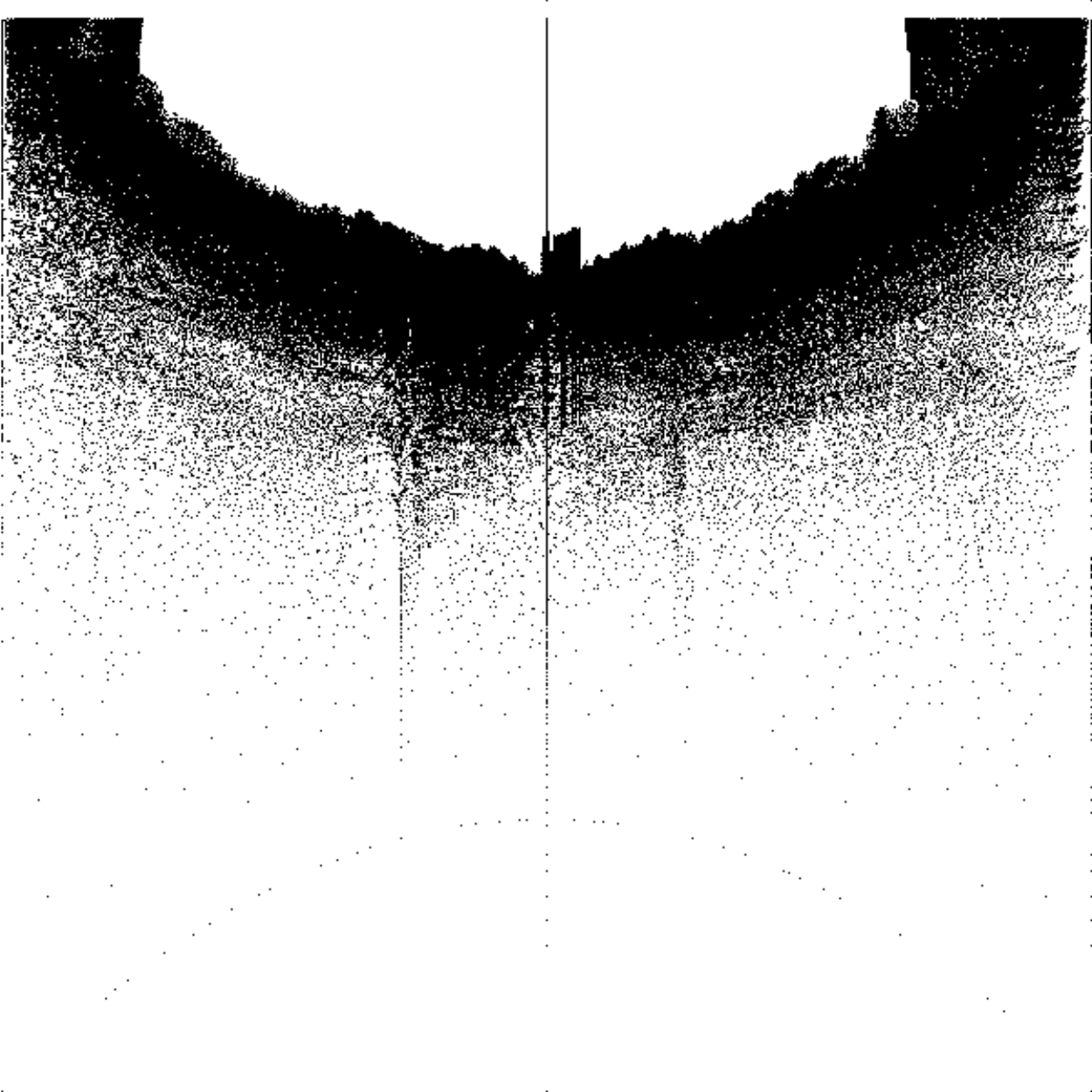}
    \caption{The teapot for Example \ref{tree}, with the unit cylinder, plotting the half where the imaginary part is non negative.}
    \label{fig:thurteapotalt}
\end{center}
\end{figure}

\begin{example}\label{tree}
    This example is inspired by an Madison Experimental Math project I supervised, with Robert Argus, Beining Mu, Anvit Thekkatte, and Zihan Zhao. Consider the following family of continuous maps on a tree: Let $T$ be a metric graph which is also a tree with 3 branches, two of which denoted as $I_0$, $I_1$, the third divided into two segments $I_2$ and $I_3$ such that $I_2$ has the center of the tree as an end point. Consider the family of continuous maps that fixes the center, send $I_0$ bijectively to $I_2\cup I_3$, $I_2$ to $I_1$, $I_1$ to $I_0$, and $I_3$ to a sub-segment of $I_0\cup I_1$. 

\begin{figure}[H]
\begin{center}
    \begin{tikzpicture}
        \draw (0, 0)--(0, -2);
        \node at (0.3, -1) {$I_0$};
        \draw (0, 0)--(-1, 2);
        \draw (0, 0)--(1, 3);
        \node at (0.5,1.5) [circle,fill,inner sep=1.5pt]{};
        \node at (-1, 1) {$I_1$};
        \node at (0.5, 0.7){$I_2$};
        \node at (1, 2){$I_3$};
    \end{tikzpicture}
    \caption{The tree in Example \ref{tree}}
    \label{fig:enter-label}
\end{center}
\end{figure}
    
    Such map corresponds to the directed graph and sign function as follows:
    \[\Gamma=\{\{0, 1, 2, 3\}, \{(0, 2), (0, 3), (1, 0), (2, 1), (3, 0), (3, 1)\}\]
    \[O(0)=O(1)=O(3)=-1, O(2)=1\]
    We can also show that this $(\Gamma, O)$ is tunable: $\wmin=3210$, and if $w$ is extremal but not admissible. Because $w^\infty\in \Sigma$, its suffix can only be $210$ or $310$ or $30$. If an extremal word ends with $30$ it has to be $(30)^n$ hence admissible. If it ends with $210$ or $310$, $\prev(w)$ is $w$ with the other $3$-suffix. Now apply Remark \ref{fliplastdig}. Figure \ref{fig:thurteapotalt} is a plot of the ``teapot'' for this example, one can see that the ``persistence'' behavior is also illustrated in the plot.
\end{example}

A question that remains is the following:

\begin{question}
Can the assumptions for tunable and unique maximal tuning be further relaxed? In particular, can Theorem \ref{dendom} be strengthened to made sufficient to cover persistent results for all veins of the Mandelbrot set, thereby strengthening the main result in \cite{lindsey2021master}.
\end{question}


\bibliographystyle{amsplain}
\bibliography{refs}

\end{document}